\documentclass[12 pt,a4]{article}
\usepackage{amssymb,amsmath,graphicx,enumerate}
\usepackage{longtable,mathrsfs}
\usepackage{hyperref}
\graphicspath {G:\Dhaigude}

\DeclareMathAlphabet{\mathpzc}{OT1}{pzc}{m}{it}

\numberwithin{equation}{section}

\newtheorem{theorem}{Theorem}[section]
\newtheorem{lemma}{Lemma}[section]

\newtheorem{definition}{Definition}[section]

\title{ Sumudu Transform of Hilfer-Prabhakar Fractional Derivatives and Applications}
\author{\bf S. K. Panchal, \,\,\,     Amol D. Khandagale,\\
\bf Pravinkumar V. Dole.\\ Department of mathematics, \\ Dr. Babasaheb Ambedkar Marathwada University,\\ Aurangabad-431004 (M.S.) India.\\
E-mail ID - drpanchalsk@gmail.com\\ kamoldsk@gmail.com\\ pvasudeo.dole@gmail.com  }
\date{\it July 2016.}
\begin{document}
\maketitle{\bf{Abstract}}
        In this paper the Sumudu transforms of Hilfer-Prabhakar fractional derivative and regularized version of Hilfer-Prabhakar fractional derivative are obtained. These results are used to obtain relation between them involving Mittag-Leffler function. Also these results are applied to solve some problems in physics. Here the solutions of problems involving Hilfer-Prabhakar fractional derivative and regularized version of Hilfer-Prabhakar fractional derivative are obtained by using Fourier and Sumudu transform techniques.\\
        
\textbf{Keywords:} 
\textit{Fractional calculus, Hilfer-Prabhakar fractional derivative, Mittag-Leffler function, Integral transforms.}\\

\textbf{2010 AMS Subject Classification:}\quad $26A33$, $42A38$, $42B10$.
\section{Introduction}
	The concept of integral transforms is originated from the Fourier integral formula. The importance of integral transforms is that they provide powerful operational method for solving initial and boundary value problems. The operational calculus of integral transform is used to solve the differential and integral equations arising in applied mathematics, mathematical physics and Engineering science. K. S. Miller, B. Ross, Podlubny Igor, A. A. Kilbas, J. J. trujillo used Laplace transform approach to solve Cauchy type fractional differential equations and F. Mainardi solved the viscoelastic problems by Laplace techniques.  Recently, Watugala G. K. \cite{WGK} introduced a new integral transform in 1990 known as Sumudu transform.  Watugala G. K. solved fractional differential equations by using Sumudu transform techniques.
\paragraph{}The Prabhakar integral \cite{PTR} is defined by modifying Riemann-Liouville integral operator by extending its kernel with a three-parameter Mittag-Leffler function. The Hilfer-Prabhakar fractional derivative and its Caputo like regularized counterpart were first introduced in \cite{RRFZ}. In this paper the Sumudu transforms of Hilfer-Prabhakar fractional derivative and regularized version of Hilfer-Prabhakar fractional derivative are obtained. These results are used to obtain relation between them involving Mittag-Leffler function and it is also used to obtain the solutions of non-homogeneous Cauchy type fractional differential equations \cite{RRFZ} in which Hilfer-Prabhakar fractional derivative and regularized version of Hilfer-Prabhakar fractional derivative are involved.
\section{Preliminaries}
\paragraph{} In this section we gives some definitions, theorems and lemma, which are used in the paper. 

\begin{definition}\cite{WGK}
Consider a set A defined as,\\
\begin{align}
A=\bigg{\lbrace}f(t)/ \exists\,\, M, \,\, \tau_{1}, \tau_{2} > 0, |f(t)|\leq Me^{\frac{|t|}{\tau_{j}}}\,\, if \,\, t\,\,\in \,\,(-1)^{j}\times[0, \infty)\bigg{\rbrace}\label{2.1}
\end{align}
For all real $t\geq 0$ the Sumudu transform of function $f(t)\in A$ is defined as,
\begin{align}
S[f(t)](u) = \int_{0}^{\infty}\frac{1}{u}e^{-\frac{t}{u}}f(t) dt,\,\,\, u\in(-\tau_{1}, \tau_{2})\label{2.2}
\end{align}
and is denoted by $F(u) = S[f(t)](u)$.
\end{definition}

\begin{definition}\cite{WGK}
The function $f(t)$ in \eqref{2.1} is called inverse Sumudu transform of $F(u)$ and is denoted by, 
\begin{align}
f(t) = S^{-1}[F(u)](t)\label{2.3}
\end{align} 
and the inversion formula for Sumudu transform is given by \cite{WGK} ,
\begin{align}
f(t) = S^{-1}[F(u)](t)= \frac{1}{2\pi i}\int_{\gamma-i\infty}^{\gamma-i\infty} \frac{1}{u} e^{\frac{t}{u}}f(u) du \label{2.4}
\end{align}
For $Re(\frac{1}{u}) > \gamma$ and $\gamma\,\,\in \,\,\mathbb{C}$.
\end{definition}

\begin{definition}\cite{PTR}
The three parameter Mittage-Leffter function introduced by Prabhakar is of the form
\begin{equation}
E_{\alpha, \beta}^{\gamma}(z) = \sum_{k=0}^{\infty}\frac{\Gamma(r+k)}{\Gamma(r) \Gamma(\alpha k+\beta)} \frac{z^{k}}{k!}\label{2.5}
\end{equation}
for $\alpha, \beta, \gamma \,\,\in \,\,\mathbb{C}$ and  $Re(\alpha) > 0$.
\end{definition}

\begin{definition}\cite{DB}
Let $f(x)$ be a function defined on $(-\infty, \infty)$ and be piecewise continuous in each finite partial interval and absolutely integrable in $(-\infty, \infty)$ then 
\begin{align}
F[f(x)](p) = \int_{-\infty}^{\infty} e^{-ipx} f(x) dx\label{2.6}
\end{align}
is called Fourier transform of $f(x)$ and is denoted by $F[f(x)](p)= f(p)$.\\
The function $f(x)$ called the inverse Fourier transform of $f(p)$, is defined as
\begin{align}
f(x) = F^{-1}[f(p)](x)= \frac{1}{2\pi}\int_{-\infty}^{\infty} f(p) e^{ipx} dp\label{2.7}
\end{align}
and is denoted by $f(x) = F^{-1}[f(p)](x)$.
\end{definition}

\begin{definition}\cite{PTR}{\bf(Prabhakar Integral)}\\
Let $f\epsilon L^{1} [0, b]; 0 < t < b < \infty$. The prabhakar integral is defined as,
\begin{align}
{\bf E}_{\rho, \mu, \omega, 0^{+}}^{\gamma} f(t) &= \int_{0}^{t} (t-y)^{\mu-1} E_{\rho, \mu}^{\gamma}[\omega(t-y)^{\rho}]f(y) dy \nonumber \\ 
&= (f*e_{\rho, \mu, \omega}^{\gamma})(t)\label{2.8}
\end{align} 
where * denote the convolution operation;\,\, $\rho, \mu, \omega, \gamma\,\,\in \,\,\mathbb{C};\,\, Re(\rho), Re(\mu) > 0$ and
\begin{align}
e_{\rho, \mu, \omega}^{\gamma}(t)= t^{\mu-1} E_{\rho, \mu}^{\gamma}[\omega t^{\rho}].\label{2.9}
\end{align}
\end{definition}

For $n\in \mathbb{N}$, we denote by $AC^{n}[a, b]$ the space of the real valued function $f(t)$ which have continuous derivative up to order $(n-1)$ on $[a, b]$ such that $f^{(n-1)}(t)$ belongs to the space of absolutely continuous functions $AC[a, b]$,\\
\begin{center}
$AC^{n}[a, b] = \bigg{\lbrace} f: [a, b]\rightarrow \mathbb{R} ; \frac{d^{n-1}}{dx^{n-1}}f(x) \in AC[a, b] \bigg{\rbrace}$.
\end{center}

\begin{definition}\cite{RRFZ} {\bf( Hilfer-Prabhakar Fractional Derivative)}\\
Let $\mu\in(0, 1);\nu\in [0, 1]$ and let $f \in L^{1}[0, b]; 0 < t < b < \infty$;
$(f*e_{\rho,(1-\nu)(1-\mu),\omega}^{-\gamma (1-\nu)})(t) \in AC^{1}[0, b]$.
The Hilfer-Prabhakar fractional derivative of $f(t)$ of order $\mu$ denoted by ${\bf D}_{\rho, \omega, 0^{+}}^{\gamma, \mu, \nu}f(t)$ is defined as,
\begin{align}
{\bf D}_{\rho, \omega, 0^{+}}^{\gamma, \mu, \nu} f(t)= \bigg{(}{\bf E}_{\rho, \nu(1-\mu), \omega, 0^{+}}^{-\gamma\nu}\frac{d}{dt}\big{(}{\bf E}_{\rho,(1-\nu)(1-\mu), \omega, 0^{+}}^{-\gamma(1-\nu)}f\big{)}\bigg{)}(t)  \label{2.10}
\end{align}
where $\gamma, \omega \in \mathbb{R}, \rho > 0$  and   ${\bf E}_{\rho, 0, \omega, 0^{+}}^{0} f = f$.
\end{definition}

In order to consider the Cauchy problems in which initial condition depending only on the function and its integer-order derivative, we use the regularized version of Hilfer-Prabhakar fractional derivative defined as below.

\begin{definition}\cite{RRFZ} {\bf(Regularized Version of Hilfer-Prabhakar Fractional \\Derivative)}\\
For $f\in AC^{1}[0, b], 0 < t < b < \infty$; $\mu\in(0, 1);\nu\in [0, 1]$; $\gamma, \omega \in \mathbb{R}, \rho > 0$. The regularized version of Hilfer-Prabhakar fractional derivative of $f(t)$ denoted by ${^C\textbf{D}}_{\rho,\omega,0^{+}}^{\gamma,\mu}f(t)$  is defined as,
\begin{align}
{^C\textbf{D}}_{\rho,\omega,0^{+}}^{\gamma,\mu}&=\bigg{(}{\bf E}_{\rho, \nu(1-\mu), \omega, 0^{+}}^{-\gamma\nu}{\bf E}_{\rho,(1-\nu)(1-\mu), \omega, 0^{+}}^{-\gamma(1-\nu)}\frac{d}{dt}f\bigg{)}(t)\label{2.11}\\
&=\bigg{(}{\bf E}_{\rho, \nu(1-\mu), \omega, 0^{+}}^{-\gamma}\frac{d}{dt}f\bigg{)}(t).\label{2.12}
\end{align}
\end{definition} 
\begin{theorem}\cite{WGK} Let $F(u)$ and $G(u)$ be Sumudu transforms of $f(t)$ and $g(t)$ respectively. The Sumudu transforms of convolution of $f$ and $g$ is
\begin{align}
S\big{[}(f*g)(t)\big{]}(u)=uF(u)G(u),\label{2.13}
\end{align}
where $(f*g)(t)=\int_{0}^{t}f(t)g(t-\tau)d\tau$. 
\end{theorem}

\begin{lemma}\cite{CDB} Let $\alpha, \beta, \lambda \in \mathbb{R}$ and $\alpha > 0, \beta > 0, n \in \mathbb{N}$. The Sumudu transform of function $e_{\rho, \mu, \omega}^{\gamma}(t)$ defined in $(2.10)$ is,
\begin{align}
S\big{[}e_{\rho, \mu, \omega}^{\gamma}(t)\big{]}(u)= \frac{u^{(\beta-1)}}{(1-\lambda u^{\alpha})^{\delta}}.\label{2.14}
\end{align}
\end{lemma}

\section{Sumudu Transform of Hilfer-Prabhakar\\ Fractional Derivatives}
\begin{lemma}
The Sumudu transform of Hilfer-Prabhakar fractional derivative \eqref{2.10} is,
\begin{align}
&S\bigg{(}{\bf D}_{\rho, \omega, 0^{+}}^{\gamma, \mu, \nu} f(t)\bigg{)}(u)=S\bigg{(}{\bf E}_{\rho, \nu(1-\mu), \omega, 0^{+}}^{-\gamma\nu}\frac{d}{dt}\big{(}{\bf E}_{\rho,(1-\nu)(1-\mu), \omega, 0^{+}}^{-\gamma(1-\nu)}f\big{)}\bigg{)}(u) \nonumber \\
&= u^{-\mu}(1-\omega u^{\rho})^{\gamma}F(u)-u^{\nu(1-\mu)-1}(1-\omega u^{\rho})^{\gamma\nu}
\bigg{[}{\bf E}_{\rho,(1-\nu)(1-\mu), \omega, 0^{+}}^{-\gamma(1-\nu)}f(t)\bigg{]}_{t=0^{+}}.\label{3.1}
\end{align}
\end{lemma}
\textbf{Proof:}
Taking Sumudu transform of Hilfer-Prabhakar fractional derivative \eqref{2.10} and using \eqref{2.8}, \eqref{2.9}, \eqref{2.13} and \eqref{2.14},
we have,
\begin{align*}
S\bigg{(}{\bf D}_{\rho,\omega,0^{+}}^{\gamma,\mu,\nu} f(t)\bigg{)}(u)&=S\bigg{(}\bigg{(}{\bf E}_{\rho,\nu(1-\mu),\omega,0^{+}}^{-\gamma\nu}\frac{d}{dt}\big{(}{\bf E}_{\rho,(1-\nu)(1-\mu),\omega,0^{+}}^{-\gamma(1-\nu)}f\big{)}\bigg{)}(t)\bigg{)}(u)\\
&=S\bigg{(} \bigg{(}e_{\rho,\nu(1-\mu),\omega}^{-\gamma\nu}*\frac{d}{dt}\big{(}{\bf E}_{\rho,(1-\nu)(1-\mu),\omega,0^{+}}^{-\gamma(1-\nu)}f\big{)}\bigg{)}(t)\bigg{)}(u)\\
&=uS\bigg{(}t^{\nu (1-\mu)-1} E_{\rho, \nu(1-\mu)}^{-\gamma\nu}(\omega t^{\rho})\bigg{)}(u)S\bigg{(}\frac{d}{dt}\big{(}{\bf E}_{\rho,(1-\nu)(1-\mu), \omega, 0^{+}}^{-\gamma(1-\nu)}f\big{)}(t)\bigg{)}(u)
\end{align*}
\begin{align*}
&=uu^{\nu(1-\mu)-1}(1-\omega u^{\rho})^{\gamma\nu}\bigg{\lbrace}\frac{S\big{[}{\bf E}_{\rho,(1-\nu)(1-\mu), \omega, 0^{+}}^{-\gamma(1-\nu)}f(t)\big{]}(u)-\big{[}{\bf E}_{\rho,(1-\nu)(1-\mu), \omega, 0^{+}}^{-\gamma(1-\nu)}f(t)\big{]}_{t=0^{+}}}{u}\bigg{\rbrace}\\
&=u^{\nu(1-\mu)-1}(1-\omega u^{\rho})^{\gamma\nu}S\bigg{(} \big{(}e_{\rho,(1-\nu)(1-\mu), \omega}^{-\gamma(1-\nu)}*f\big{)}(t)\bigg{)}(u)\\
&\qquad-u^{\nu(1-\mu)-1}(1-\omega u^{\rho})^{\gamma\nu}\bigg{[}{\bf E}_{\rho,(1-\nu)(1-\mu), \omega, 0^{+}}^{-\gamma(1-\nu)}f(t)\bigg{]}_{t=0^{+}}\\
&=u^{\nu(1-\mu)-1}(1-\omega u^{\rho})^{\gamma\nu}uS\bigg{(}t^{(1-\nu)(1-\mu)-1} E_{\rho,(1-\nu)(1-\mu)}(\omega t^{\rho})\bigg{)}(u) S[f(t)](u)\\
&\qquad-u^{\nu(1-\mu)-1}(1-\omega u^{\rho})^{\gamma\nu}\bigg{[}{\bf E}_{\rho,(1-\nu)(1-\mu),\omega, 0^{+}}^{-\gamma(1-\nu)}f(t)\bigg{]}_{t=0^{+}}\\
&= u^{\gamma(1-\mu)}(1-\omega u^{\rho})^{\gamma\nu}u^{(1-\nu)(1-\mu)-1}(1-\omega u^{\rho})^{\gamma(1-\nu)}S[f(t)](u)\\
&\qquad-u^{\nu(1-\mu)-1}(1-\omega u^{\rho})^{\gamma\nu}\bigg{[}{\bf E}_{\rho,(1-\nu)(1-\mu),\omega, 0^{+}}^{-\gamma(1-\nu)}f(t)\bigg{]}_{t=0^{+}}.
\end{align*}
Thus we get the required result \eqref{3.1}.

\begin{lemma}\label{lem3.2}
The Sumudu transforms of the regularized version of Hilfer-Prabhakar fractional derivative \eqref{2.11} of order $\mu$ is,
\begin{align}
S\bigg{(}{^C\textbf{D}}_{\rho, \omega, 0^{+}}^{\gamma, \mu} f(t)\bigg{)}(u)&=u^{-\mu}(1-\omega u^{\rho})^{\gamma}\bigg{(}F(u)-f(0^{+})\bigg{)} \nonumber\\
&=u^{-\mu}(1-\omega u^{\rho})^{\gamma}F(u)-u^{-\mu}(1-\omega u^{\rho})^{\gamma}f(0^{+}).\label{3.2}
\end{align}
\end{lemma}
\textbf{Proof:} Taking Sumudu transforms of regularized version of Hilfer-Prabhakar fractional derivative \eqref{2.11} of order $\mu$ and using \eqref{2.8}, \eqref{2.9}, \eqref{2.13} and \eqref{2.14}. We have,
\begin{align*}
&S\bigg{(}{^C\textbf{D}}_{\rho, \omega, 0^{+}}^{\gamma, \mu} f(t)\bigg{)}(u)=S\bigg{(}\bigg{(}{\bf E}_{\rho,\nu(1-\mu),\omega,0^{+}}^{-\gamma\nu}\big{(}{\bf E}_{\rho,(1-\nu)(1-\mu),\omega,0^{+}}^{-\gamma(1-\nu)}\frac{d}{dt}f\big{)}\bigg{)}(t)\bigg{)}(u)\\
&=S\bigg{(} \bigg{(}e_{\rho,\nu(1-\mu),\omega}^{-\gamma\nu}*\big{(}{\bf E}_{\rho,(1-\nu)(1-\mu),\omega,0^{+}}^{-\gamma(1-\nu)}\frac{d}{dt}f\big{)}\bigg{)}(t)\bigg{)}(u)\\
&=uS\bigg{(}t^{\nu (1-\mu)-1} E_{\rho, \nu(1-\mu)}^{-\gamma\nu}(\omega t^{\rho})\bigg{)}(u)S\bigg{(}{\bf E}_{\rho,(1-\nu)(1-\mu),\omega,0^{+}}^{-\gamma(1-\nu)}\frac{d}{dt}f(t)\bigg{)}(u)\\
&=uu^{\nu(1-\mu)-1}(1-\omega u^{\rho})^{\gamma\nu}S\bigg{(}\bigg{(}e_{\rho,(1-\nu)(1-\mu),\omega}^{-\gamma(1-\nu)}*\frac{d}{dt}f\bigg{)}(t)\bigg{)}(u)\\
&=u^{\nu(1-\mu)}(1-\omega u^{\rho})^{\gamma\nu}uS\bigg{(}t^{(1-\nu)(1-\mu)-1}E_{\rho,(1-\nu)(1-\mu)}^{-\gamma(1-\nu)}(\omega t^{\rho})\bigg{)}(u)S\bigg{(}\frac{d}{dt}f(t)\bigg{)}(u)\\
&=u^{\nu(1-\mu)}(1-\omega u^{\rho})^{\gamma\nu}uu^{(1-\nu)(1-\mu)-1}(1-\omega u^{\rho})^{\gamma(1-\nu)}\bigg{\lbrace}\frac{S\big{[}f(t)\big{]}(u)-f(0^{+})}{u}\bigg{\rbrace}
\end{align*}
Thus we get the required result \eqref{3.2}.\\

\textbf{Alternating Proof of Lemma \eqref{lem3.2}:} Taking Sumudu transforms of regularized version of Hilfer-Prabhakar fractional derivative  \eqref{2.12} of order $\mu$ and using \eqref{2.8}, \eqref{2.9}, \eqref{2.13}, \eqref{2.14}. We have,
\begin{align*}
S\bigg{(}{^C\textbf{D}}_{\rho, \omega, 0^{+}}^{\gamma, \mu} f(t)\bigg{)}(u)&=S\bigg{(}{\bf E}_{\rho, \nu(1-\mu), \omega, 0^{+}}^{-\gamma}\frac{d}{dt}f\bigg{)}(u)\\
&=S\bigg{(}\bigg{(}e_{\rho,(1-\mu),\omega}^{-\gamma}*\frac{d}{dt}f\bigg{)}(t)\bigg{)}(u)\\
&=uS\bigg{(}t^{-\mu}E_{\rho,(1-\mu)}^{-\gamma}(\omega t^{\rho})\bigg{)}(u)S\bigg{(}\frac{d}{dt}f(t)\bigg{)}(u)\\
&=uu^{-\mu}(1-\omega u^{\rho})^{\gamma}\bigg{\lbrace}\frac{S\big{[}f(t)\big{]}(u)-f(0^{+})}{u}\bigg{\rbrace}\\
&=u^{-\mu}(1-\omega u^{\rho})^{\gamma}\bigg{(}F(u)-f(0^{+})\bigg{)}.
\end{align*}

Observing that for absolutely continuous function $f\in AC^{1}[0, b]$,
\begin{align}
\big{[}{\bf E}_{\rho,(1-\nu)(1-\mu),\omega, 0^{+}}^{-\gamma(1-\nu)}f(t)\big{]}_{t=0^{+}}=0 \label{3.3}
\end{align} 
then the result \eqref{3.1} becomes,
\begin{align}
S\bigg{(}{\bf{D}}_{\rho, \omega, 0^{+}}^{\gamma, \mu, \nu} f(t)\bigg{)}(u)=u^{-\mu}(1-\omega u^{\rho})^{\gamma}F(u)\label{3.4}
\end{align} 
Substituting this value \eqref{3.4} in \eqref{3.2}, we get,
\begin{align}
S\bigg{(} {^CD}_{\rho, \omega, 0^{+}}^{\gamma, \mu} f(t)\bigg{)}(u)=S\bigg{(}{\bf{D}}_{\rho, \omega, 0^{+}}^{\gamma, \mu, \nu} f(t)\bigg{)}(u)-u^{-\mu}(1-\omega u^{\rho})^{\gamma}f(0^{+})\label{3.5}
\end{align}
taking inverse Sumudu transform of \eqref{3.5}, we get the relation between Hilfer-Prabhakar fractional derivative and regularized version of Hilfer-Prabhakar fractional derivative in terms of Mittag-leffter function as below,
\begin{align}
{^CD}_{\rho, \omega, 0^{+}}^{\gamma, \mu} f(t)={\bf D}_{\rho, \omega, 0^{+}}^{\gamma, \mu, \nu} f(t)-t^{-\mu}E_{\rho,(1-\mu)}^{-\gamma}(\omega t^{\rho})f(0^{+}),\label{3.6}
\end{align}
for $f\in AC^{1}[0, b]$.

\section{Applications of Hilfer-Prabhakar Fractional\\ Derivatives using Sumudu transform and Fourier transform techniques}
\begin{theorem}
The solution of Cauchy problem
\begin{align}
{\bf D}_{\rho, \omega, 0^{+}}^{\gamma, \mu, \nu} y(x)&=\lambda {\textbf{E}}_{\rho,\mu, \omega,0^{+}}^{\delta}y(x)+f(x),\label{4.1}\\
\big{[}{\bf E}_{\rho,(1-\nu)(1-\mu),\omega, 0^{+}}^{-\gamma(1-\nu)}f(t)\big{]}_{t=0^{+}}&=K, \label{4.2}
\end{align}
where $f(x)\in L^{1}[0,\infty)$; $\mu\in(0, 1)$, $\nu\in[0, 1]$; $\omega, \lambda\in \mathbb{C}$; $x, \rho > 0$, $K, \gamma, \delta\geq 0$, is
\begin{align}
y(x)=K\sum_{n=0}^{\infty}\lambda^{n}x^{\nu(1-\mu)+\mu(2n+1)-1}E_{\rho,\nu(1-\mu)+\mu(2n+1)}^{\gamma-\gamma\nu+n(\delta+\gamma)}(\omega x^{\rho})+\sum_{n=0}^{\infty}{\textbf{E}}_{\rho,(2n+1)\mu,\omega,0^{+}}^{\gamma+n(\delta+\gamma)}f(x).\label{4.3}
\end{align}
\end{theorem}
\textbf{Proof:}
Let $Y(u)$ and $F(u)$ denote the Sumudu transform of $ y(x) $ and $ f(x) $ respectively, Now taking Sumudu transform of \eqref{4.1} and using \eqref{2.8}, \eqref{2.9}, \eqref{2.13}, we get
\begin{align*}
S\bigg{(}{\bf D}_{\rho,\omega,0^{+}}^{\gamma,\mu,\nu} y(x)\bigg{)}(u)&=\lambda S\bigg{(}{\textbf{E}}_{\rho,\mu,\omega,0^{+}}^{\delta}y(x)\bigg{)}(u)+S\bigg{(}f(x)\bigg{)}(u) \\
&=\lambda S\bigg{(}\big{(}e_{\rho,\mu,\omega}^{\delta}*y\big{)}(x)\bigg{)}(u)+F(u)\\
&=\lambda u S\bigg{(}x^{\mu-1}E_{\rho,\mu}^{\delta}(\omega x^{\rho})\bigg{)}(u)S\big{(} y(x)\big{)}(u)+F(u)
\end{align*}
From \eqref{2.14},\eqref{3.1} and \eqref{4.2}, we get,
\begin{align*}
u^{-\mu}(1-\omega u^{\rho})^{\gamma}Y(u)-u^{\nu(1-\mu)-1}(1-\omega u^{\rho})^{\gamma\nu}K=\lambda u u^{\mu-1}(1-\omega u^{\rho})^{-\delta}Y(u)+F(u)\\
\bigg{(}u^{-\mu}(1-\omega u^{\rho})^{\gamma}-\lambda  u^{\mu}(1-\omega u^{\rho})^{-\delta}\bigg{)}Y(u)=u^{\nu(1-\mu)-1}(1-\omega u^{\rho})^{\gamma\nu}K+F(u)
\end{align*}
\begin{align*}
Y(u)&=\bigg{(}\frac{Ku^{\nu(1-\mu)-1}(1-\omega u^{\rho})^{\gamma\nu}}{u^{-\mu}(1-\omega u^{\rho})^{\gamma}}+\frac{F(u)}{u^{-\mu}(1-\omega u^{\rho})^{\gamma}}\bigg{)}\frac{1}{\bigg{(}1-\frac{\lambda  u^{\mu}(1-\omega u^{\rho})^{-\delta}}{u^{-\mu}(1-\omega u^{\rho})^{\gamma}}\bigg{)}}\\
&=\bigg{(}Ku^{\nu(1-\mu)+\mu-1}(1-\omega u^{\rho})^{\gamma\nu-\gamma}+F(u)u^{\mu}(1-\omega u^{\rho})^{-\gamma}\bigg{)}\sum_{n=0}^{\infty}\lambda^{n}u^{2n\mu}(1-\omega u^{\rho})^{-n(\delta+\gamma)}\\
&=K\sum_{n=0}^{\infty}u^{\nu(1-\mu)+\mu(2n+1)-1}(1-\omega u^{\rho})^{-n(\delta+\gamma)+\gamma\nu-\gamma}+F(u)\sum_{n=0}^{\infty}u^{\mu(2n+1)}(1-\omega u^{\rho})^{-n(\delta+\gamma)-\gamma}.
\end{align*}
Taking inverse Sumudu transform on both side of above equation, we get the required solution \eqref{4.3}.
\begin{theorem}
The solution of Cauchy problem
\begin{align}
{^C\textbf{D}}_{\rho, -\omega,0^{+}}^{\gamma, \mu}G(v,t)&=-\lambda(1-v)G(v,t),\qquad for \mid v\mid\leq 1 \label{4.4}\\
G(v,0)&=1,\label{4.5}
\end{align}
with $t > 0$, $\phi, \lambda > 0$, $\gamma \geq 0$, $0 < \rho \leq 1$, $0 < \mu \leq 1$, is
\begin{align}
G(v,t)=\sum_{n=o}^{\infty}(-\lambda)^{n}(1-v)^{n}t^{n\mu}E_{\rho,n\mu+1}^{n\gamma}(-\omega t^{\rho}).\label{4.6}
\end{align}
\end{theorem}
\textbf{Proof:} Taking Sumudu transform of \eqref{4.4} with respect to $t$ and using \eqref{3.2} and \eqref{4.5}.
\begin{align*}
S\bigg{(}{^C\textbf{D}}_{\rho, -\omega,0^{+}}^{\gamma, \mu}G(v,t)\bigg{)}(v,q)&=-\lambda(1-v)S\bigg{(}G(v,t)\bigg{)}(v,q)\\
q^{-\mu}(1+\omega q^{\rho})^{\gamma}S\bigg{(}G(v,t)\bigg{)}(v,q)-q^{-\mu}(1+\omega q^{\rho})^{\gamma}G(v,0)&=-\lambda(1-v)S\bigg{(}G(v,t)\bigg{)}(v,q)\\
\bigg{[}q^{-\mu}(1+\omega q^{\rho})^{\gamma}+\lambda(1-v)\bigg{]}S\bigg{(}G(v,t)\bigg{)}(v,q)&=q^{-\mu}(1+\omega q^{\rho})^{\gamma}
\end{align*}
\begin{align*}
S\bigg{(}G(v,t)\bigg{)}(v,q)&=\frac{q^{-\mu}(1+\omega q^{\rho})^{\gamma}}{q^{-\mu}(1+\omega q^{\rho})^{\gamma}\bigg{(}1+\frac{\lambda(1-v)}{q^{-\mu}(1+\omega q^{\rho})^{\gamma}}\bigg{)}}\\
S\bigg{(}G(v,t)\bigg{)}(v,q)&=\sum_{n=0}^{\infty}(-\lambda)^{n}(1-v)^{n}q^{n\mu}(1+\omega q^{\rho})^{-n\gamma}
\end{align*}
Taking inverse Sumudu transform on both side of above equation, we get the required solution \eqref{4.6}.
\begin{theorem}
The solution of Cauchy problem
\begin{align}
{\bf D}_{\rho, \omega, 0^{+}}^{\gamma, \mu, \nu} u(x,t)&=K\frac{\partial^{2}}{\partial x^{2}}u(x,t),\label{4.7}\\
\bigg{[}{\bf E}_{\rho,(1-\nu),(1-\mu), \omega, 0^{+}}^{-\gamma(1-\nu)}u(x,t)\bigg{]}_{t=0^{+}}&=g(x),\label{4.8}\\
\lim_{x\to\pm\infty} u(x,t)&=0,\label{4.9}
\end{align}
with $\mu\in(0,1)$, $\nu\in[0, 1]$; $x, \omega\in \mathbb{R}$; $t, \rho > 0$; $K, \gamma \geq 0$, is
\begin{align}
u(x,t)=\frac{1}{2\pi}\int_{-\infty}^{\infty}dp\, e^{ipx} \widehat{g}(p)\sum_{n=0}^{\infty}(-k)^{n}p^{2n}t^{\mu(n+1)-\nu(\mu-1)-1}E_{\rho,\mu(n+1)-\nu(\mu-1)}^{\gamma(n+1-\nu)}(\omega t^{\rho}).\label{4.10}
\end{align}
\end{theorem}
\textbf{Proof:} Let $\overline{u}(x, q)$ denote the Sumudu transform of $u(x, t)$ and $\widehat{u}(p, t)$ denote the Fourier transform of $u(x, t)$. Taking Fourier-Sumudu transform of \eqref{4.7} and using \eqref{3.1}, \eqref{4.8}, \eqref{4.9} we get,
\begin{align*}
q^{-\mu}(1-\omega q^{\rho})^{\gamma}\quad\widehat{\overline{u}}(p, q)-q^{\nu(1-\mu)-1}(1-\omega q^{\rho})^{\gamma\nu}\widehat{g}(p)&=-K p^{2}\quad\widehat{\overline{u}}(p,q)\\
\bigg{(}q^{-\mu}(1-\omega q^{\rho})^{\gamma}+K p^{2}\bigg{)}\quad\widehat{\overline{u}}(p,q)&=q^{\nu(1-\mu)-1}(1-\omega q^{\rho})^{\gamma\nu}\widehat{g}(p)
\end{align*}
\begin{align}
\widehat{\overline{u}}(p,q)&=\frac{q^{\nu(1-\mu)-1}(1-\omega q^{\rho})^{\gamma\nu}}{q^{-\mu}(1-\omega q^{\rho})^{\gamma}\bigg{(}1+\frac{K p^{2}}{q^{-\mu}(1-\omega q^{\rho})^{\gamma}} \bigg{)}}\widehat{g}(p) \nonumber\\
\widehat{\overline{u}}(p,q)&=\widehat{g}(p)q^{\mu+\nu(1-\mu)-1}(1-\omega q^{\rho})^{\gamma\nu-\gamma}\sum_{n=0}^{\infty}(-K)^{n} p^{2n}q^{n\mu}(1-\omega q^{\rho})^{-n\gamma} \nonumber\\
\widehat{\overline{u}}(p,q)&=\widehat{g}(p)\sum_{n=0}^{\infty}(-K)^{n} p^{2n}q^{\nu(1-\mu)+\mu(n+1)-1}(1-\omega q^{\rho})^{\gamma\nu-\gamma-n\gamma}\label{4.11}
\end{align}
Taking inverse Fourier transform of above equation \eqref{4.11} we get,
\begin{align}
\overline{u}(x,q)=\frac{1}{2\pi}\int_{-\infty}^{\infty}dp\quad e^{ipx}\widehat{g}(p)\sum_{n=0}^{\infty}(-K)^{n} p^{2n}q^{\nu(1-\mu)+\mu(n+1)-1}(1-\omega q^{\rho})^{\gamma\nu-\gamma-n\gamma}\label{4.12}
\end{align}
Taking inverse Sumudu transform of above equation \eqref{4.12}, we get the required solution \eqref{4.10}.
\begin{theorem}
The solution of Cauchy problem
\begin{align}
{^C\textbf{D}}_{\rho, \omega, 0^{+}}^{\gamma, \mu} u(x,t)&=K\frac{\partial^{2}}{\partial x^{2}}u(x,t),\label{4.13}\\
u(x,0^{+})&=g(x),\label{4.14}\\
\lim_{x\to\pm\infty} u(x,t)&=0,\label{4.15}
\end{align}
with $\mu\in(0,1)$, $x, \omega\in \mathbb{R}$, $t, K, \rho > 0$, $\gamma \geq 0$, is
\begin{align}
u(x,t)=\frac{1}{2\pi}\int_{-\infty}^{\infty}dp\, e^{ipx} \widehat{g}(p)\sum_{n=0}^{\infty}(-k)^{n}p^{2n}t^{n\mu}E_{\rho,n\mu+1}^{n\gamma}(\omega t^{\rho}).\label{4.16}
\end{align}
\end{theorem}
\textbf{Proof:} Let $\overline{u}(x, q)$ denote the Sumudu transform of $u(x, t)$ and $\widehat{u}(p, t)$ denote the Fourier transform of $u(x, t)$. Now taking Fourier-Sumudu transform of \eqref{4.13} and using \eqref{3.2}, \eqref{4.14}, \eqref{4.15} we get,
\begin{align*}
q^{-\mu}(1-\omega q^{\rho})^{\gamma}\quad\widehat{\overline{u}}(p,q)-q^{-\mu}(1-\omega q^{\rho})^{\gamma}\quad\widehat{g}(p)&=-K\,p^{2}\quad\widehat{\overline{u}}(p, q)\\
\bigg{(}q^{-\mu}(1-\omega q^{\rho})^{\gamma}+K\,p^{2}\bigg{)}\widehat{\overline{u}}(p, q)&=q^{-\mu}(1-\omega q^{\rho})^{\gamma}\quad\widehat{g}(p)
\end{align*}
\begin{align}
\widehat{\overline{u}}(p, q)&=\frac{q^{-\mu}(1-\omega q^{\rho})^{\gamma}}{q^{-\mu}(1-\omega q^{\rho})^{\gamma}\bigg{(}1+\frac{K\,p^{2}}{q^{-\mu}(1-\omega q^{\rho})^{\gamma}}\bigg{)}}\quad\widehat{g}(p) \nonumber \\
\widehat{\overline{u}}(p, q)&=\widehat{g}(p)\sum_{n=0}^{\infty}(-K)^{n} p^{2n} q^{n\mu}(1-\omega q^{\rho})^{-n\gamma}\label{4.17}
\end{align}
Taking inverse Fourier transform of \eqref{4.17} , we get,
\begin{align}
\overline{u}(x,q)=\frac{1}{2\pi}\int_{-\infty}^{\infty}dp\,e^{ipx}\widehat{g}(p)\sum_{n=0}^{\infty}(-K)^{n} p^{2n} q^{n\mu}(1-\omega q^{\rho})^{-n\gamma}\label{4.18}
\end{align}
Taking inverse Sumudu transform of \eqref{4.18}, we get required solution \eqref{4.16}.

\end{document}